\documentclass[12pt]{amsart}
\usepackage{amsfonts,amsmath,amssymb,amscd,amsthm,url,color}
\long\def\comment#1\endcomment{}

\newcommand{\red}[1]{#1}

\def\Z{{\mathbb Z}} 
\def\C{{\mathbb C}} \def\Q{{\mathbb Q}} 

    \theoremstyle{theorem}
         \newtheorem{Theorem}{Theorem}
         \newtheorem{Lemma}[Theorem]{Lemma}
         
         \newtheorem{Proposition}[Theorem]{Proposition}
    \theoremstyle{definition}
         \newtheorem{Remark}[Theorem]{Remark}

\begin{document}

\title{A short elementary proof of the insolvability of the equation of degree 5}

\author{A. Skopenkov}

\thanks{Moscow Institute of Physics and Technology and Independent University of Moscow.
E-mail: skopenko@mccme.ru .
Homepage: \url{https://users.mccme.ru/skopenko/}.
\newline
Supported in part by the D. Zimin Dynasty Foundation and Simons-IUM fellowship.
\newline
This text is based on the author's lectures at Moscow `olympic' school (2015-2019), math circle `Olympiades and Mathematics' (2015-2018) and Summer Conference of Tournament of Towns (2018).
I am grateful to J. Bewersdorff, I. Bogdanov, G. Chelnokov, A. Esterov, A. Kanunnikov, E. Kogan, F. Petrov and V. Volkov for useful discussions, and to A.B. Sossinsky for correcting English in some historical remarks.}

\date{}

\maketitle

\begin{abstract}
We present short elementary proofs of the well-known Ruffini-Abel-Galois theorems on insolvability of algebraic equations in radicals.
These proofs are obtained from existing expositions by stripping away material not required for the proofs
(but presumably required elsewhere).
In particular, we do not use the terms `Galois group' and even `group'.
However, our presentation is a good way to learn (or to recall) a starting idea of Galois theory: the symmetry of a polynomial of several variables is decreased when a radical is extracted.
So the note provides a bridge (by showing that there is no gap) between elementary mathematics and Galois theory.
The note is accessible to students familiar with polynomials, complex numbers and permutations;
so the note might be interesting easy reading for professional mathematicians.
\end{abstract}


\section*{Statement and proof}

Here we present statement and proof, see discussions in the next section, \cite[\S2]{Sk17}, \cite[\S8]{Sk20}, \cite[\S5]{ZSS}.

Take a subset $X\subset\C$ containing number $1$.
A complex number $a$ is called {\bf expressible by radicals from $X$} if $a$ can be obtained from $X$ using operations of addition, subtraction, multiplication,
division by a non-zero number and taking the $n$-th root, where $n$ is a positive integer.
Or, in other words, if some set containing $a$ can be obtained
from $X$
using the following operations.
To a given set $M\subset\C$  containing numbers $x,y\in M$  one can add
$$\text{numbers}\quad x+y,x-y,xy,\quad\text{number $x/y$ when }y\ne0,$$
$$\text{and any number $r\in\C$ such that $r^n=x$ for some integer }n>0.$$



\begin{Theorem}
\label{t:galois}
For every $n\ge5$ there are $a_0,\ldots,a_{n-1}\in \C$ such that no root of the equation $x^n+a_{n-1}x^{n-1}+\ldots+a_1x+a_0=0$ is
expressible by radicals from $\{1,a_0,\ldots,a_{n-1}\}$.
\end{Theorem}

In this note we present a short
direct proof of this result, cf. Remark \ref{r:stat}.
By a proof I mean a `real' proof using elementary mathematics, not a deduction from some non-elementary results whose proofs contain
(in a longer and less motivated form) a `real' proof.
This proof is interesting because it contains an idea of an algorithm for recognition of solvability in radicals (Remark \ref{r:alg}).



We first prove a weaker insolvability result, the Ruffini Theorem \ref{t:rufnum} below.
The idea of proof is that certain symmetry is kept through extraction of radical (Lemma \ref{l:radq} below).
We then deduce
Theorem \ref{t:galois} from the Ruffini Theorem \ref{t:rufnum}.
The deduction is based on the Rationalization Lemma \ref{l:impro} below, which is also important because it illustrates one of the main ideas of Galois' (and maybe Abel's) work: if an equation is solvable in radicals at all, then it is solvable in radicals using {\it Lagrange resolvents}, see Remark \ref{r:alg}.
The Rationalization Lemma \ref{l:impro} uses the Conjugation Lemma \ref{l:irre}.b below, which
introduces the idea of a {\it field automorphism} in the simple particular case of {\it conjugation} sufficient for
Theorem \ref{t:galois}.

Denote
$$\varepsilon_k:=\cos\dfrac{2\pi}k+i\sin\dfrac{2\pi}k,\quad \Q_\varepsilon:=\bigcup\limits_{k=3}^\infty\Q(\varepsilon_3,\varepsilon_4,\ldots,\varepsilon_k)
\quad\text{and}\quad \vec y:=(y_1,\ldots,y_n)\quad\red{\text{but}\quad \vec a:=(a_{n-1},\ldots,a_0)}.$$


We use the standard notation $F[u_1,\ldots,u_n]$ and $F(u_1,\ldots,u_n)$ for the sets of polynomials and rational functions (i.e., formal ratios of polynomials) with coefficients in $F$.
Define {\bf an extension of a field $F\subset\C$ by numbers $r_1,\ldots,r_s\in\C$} as $$F(r_1,\ldots,r_s):=\{P(r_1,\ldots,r_s)\ :\ P\in F(u_1,\ldots,u_s)\}.$$
If for every $j=1,\ldots,s$ there is an integer $k_j$ such that
$r_j^{k_j}\in \red{F_{j-1}:}=F(r_1,\ldots,r_{j-1})$,
then the extension is called {\bf a radical extension}.

\begin{Theorem}[Ruffini]\label{t:rufnum}
For every $n\ge5$
there are $a_0,\ldots,a_{n-1}\in \C$ such that no root of the equation $x^n+a_{n-1}x^{n-1}+\ldots+a_1x+a_0=0$
is contained in any radical extension of $\Q_\varepsilon(\vec a)$ contained in $\Q_\varepsilon(\vec x)$, where $x_1,\ldots,x_n$ are the roots of the equation.\footnote{In order to understand the main idea one can replace $\Q_\varepsilon(\vec x)$ by $\Q[\vec x]$.
Theorem \ref{t:rufnum} holds, with analogous proof, if we replace $\Q_\varepsilon$ by any countable field.}
\end{Theorem}

For a permutation $\alpha$ denote
$$\vec u_\alpha:=(u_{\alpha(1)},\ldots,u_{\alpha(n)}).$$
A rational function $P\in\C(\vec u)$ is {\bf even-symmetric} if $P(\vec u)=P(\vec u_{(abc)})$ for every cycle $(abc)$ of length 3.\footnote{A permutation is {\it even} if it is a composition of an even number of transpositions.
Being even-symmetric is equivalent to $P(\vec u)=P(\vec u_\alpha)$ for every even permutation $\alpha$ of $\{1,\ldots,n\}$.
Indeed, any even permutation is composition of permutations of the form $(ab)(bc)=(abc)$ and $(ab)(cd)=(abc)(bcd)$.}

\begin{Lemma}\label{l:radq}
If $P$ is a rational function of $n\ge5$ variables with coefficients in $\C$, and $P^k$ is even-symmetric for some integer $k$, then $P$ is even-symmetric.
\end{Lemma}

\begin{proof}
We may assume that $k$ is a prime and $P\ne0$.

Let $a,b,c,d,e$ be arbitrary different elements of $\{1,\ldots,n\}$.

First assume that $k\ne3$.
Since
$$P^k(\vec u)=P^k(\vec u_{(abc)}),\quad\text{we have}\quad
\prod_{j=0}^{k-1}(P(\vec u)-\varepsilon_k^j P(\vec u_{(abc)}))=0.$$
Since $P$ is a non-zero rational function, there is
$$j=j(abc)\in\Z\quad\text{such that}\quad P(\vec u)=\varepsilon_k^j P(\vec u_{(abc)}).$$
Then
$$P(\vec u)= \varepsilon_k^jP(\vec u_{(abc)})= \varepsilon_k^{2j}P(\vec u_{(abc)^2})= \varepsilon_k^{3j}P(\vec u).$$
Hence $j\equiv0\mod k$, i.e. $P(\vec u)=P(\vec u_{(abc)})$.

For $k=3$ let $\sigma:=(ab)(de)=(abe)(bed)$.
Then analogously $0\equiv j(\sigma^2)\equiv2j(\sigma)\mod3.$
Hence $j(\sigma)\equiv0\mod3$, i.e. $P(\vec u)=P(\vec u_\sigma)$.
Analogously  $P(\vec u)=P(\vec u_{(ac)(de)})$.
Since  $(ab)(de)(ac)(de)=(abc)$, we have $P(\vec u)=P(\vec u_{(abc)})$.
\end{proof}

\begin{proof}[Proof of the Ruffini Theorem \ref{t:rufnum}]
Numbers $x_1,\ldots,x_n\in \C$ are called {\it algebraically independent} over $\Q_\varepsilon$ if
$P(\vec x)\ne0$ for every non-zero polynomial $P$ with coefficients in $\Q_\varepsilon$.
By induction on $n$ {\it there are $n$ algebraically independent numbers  $x_1,\ldots,x_n$ over $\Q_\varepsilon$.}
The inductive step follows because $\C$ is uncountable, while the set of roots of polynomials with coefficients in
$\Q_\varepsilon(x_1,\ldots,x_{n-1})$ is countable.

Denote the coefficients of the unitary polynomial with roots $x_1,\ldots,x_n$ by
$$a_{n-1}:=-(x_1+\ldots+x_n),\quad \ldots,\quad a_0=(-1)^nx_1\cdot\ldots\cdot x_n.$$
Assume to the contrary that there is a radical extension $\Q_\varepsilon(\vec a,r_1,\ldots,r_s)$ of
$\Q_\varepsilon(\red{\vec a})$, which both contains $x_1$ and is contained in $\Q_\varepsilon(\vec x)$.
Using Lemma \ref{l:radq}, by induction on $j$ we obtain that $r_j$ is the value at $\vec x$ of an even-symmetric rational function for every $j=1,\ldots,s$.
Since $x_1\in \Q_\varepsilon(\vec a,r_1,\ldots,r_s)$, we see that $x_1$ is also the value at $\vec x$ of an even-symmetric rational function.
Since $x_1,\ldots,x_n$ are algebraically independent over $\Q_\varepsilon$, the only such rational function is $P(u_1,\ldots,u_n)=u_1$.
This is not even-symmetric because the cycle $(123)$ carries $u_1$ to $u_2\ne u_1$.
A contradiction.
\end{proof}


Theorem \ref{t:galois} follows by the Ruffini Theorem \ref{t:rufnum} and the Abel/Ruffini Theorem \ref{t:aberuf}.

\begin{Theorem}[Abel/Ruffini]\label{t:aberuf} \red{Take any $a_0,\ldots,a_{n-1}\in \C$
such that the equation $x^n+a_{n-1}x^{n-1}+\ldots+a_1x+a_0=0$ has $n$ distinct
roots $x_1,\ldots,x_n$.
If there is a radical extension $\Q_\varepsilon(\vec a)=F_0\subset\ldots\subset F_s$ such that $x_1\in F_s$, then}
there is a radical extension $\Q_\varepsilon(\vec a)=Q_0\subset\ldots\subset Q_t$ such that
$x_1\in Q_t\subset \Q_\varepsilon(\vec x)$.\footnote{Here `$x_1\in Q_t\subset \Q_\varepsilon(\vec x)$' can be replaced by `$Q_t=F_s\cap\Q_\varepsilon(\vec x)$'. The proof is analogous.}
\end{Theorem}

Proof of Theorem \ref{t:aberuf} does not use permutations.
We construct $Q_1,\ldots,Q_t$ inductively using the following lemma, which asserts that if $F[r]$ contains more (values of) rational functions of $x_1,\ldots,x_n$ with coefficients in $\Q_\varepsilon$ than $F$, then we may assume that $r$ itself is such an `excess' rational function.


\begin{Lemma}[Rationalization]\label{l:impro} Let $n$ be an integer, $x_1,\ldots,x_n,r\in\C$ numbers,
$k$ a prime and $F\subset\C$ a field containing elementary symmetric polynomials of $x_1,\ldots,x_n$ and also   $\varepsilon_k$, $r^k$ but not $r$.
If $F(r)\cap\Q_\varepsilon(\vec x)\not\subset F$, then there is $\rho\in\Q_\varepsilon(\vec x)$ such that $\rho^k\in F$ and $F(\rho)=F(r)$.
\end{Lemma}

For a proof we need two lemmas.

\begin{Lemma}\label{l:irre}
Let $k$ be a prime, $r\in\C$ a number and $F\subset\C$ a field containing $\varepsilon_k$, $r^k$ but not $r$.

(a) (Irreducibility) Then the polynomial $z^k-r^k\in F[z]$ is irreducible over $F$.

(b) (Conjugation) If $Q\in F[z]$ a polynomial and $Q(r)=0$, then $Q(r\varepsilon_k^j)=0$ for every $j=1,\ldots,k-1$.
\end{Lemma}

\begin{proof}[Proof of (a)]
All the roots of the polynomial $z^k- r^k$ are $r,r\varepsilon_k,r\varepsilon_k^2,\ldots,r\varepsilon_k^{k-1}$.
Then the free coefficient of a factor of $z^k-r^k$ is the product of some $m$ of these roots.
Since $\varepsilon_k$, we obtain $r^m\in F$.
For a proper factor, if it existed, $0<m<k$.
Since $k$ is a prime, $ka+mb=1$ for some integers $a,b$.
Then $r=(r^k)^a(r^m)^b\in F$. A contradiction.
\end{proof}

\begin{proof}[Proof of (b)]
Since $Q(r)=0$, the remainder of division of $Q$ by $z^k-r^k$ assumes value 0 at $r$.
Since the degree of this remainder is less than $k$, by (a) this remainder is zero.
Thus $Q$ is divisible by $z^k-r^k$.
For every $j=0,1,\ldots,k-1$ since $(r\varepsilon_k^j)^k= r^k$, we obtain $Q(r\varepsilon_k^j)=0$.
\end{proof}

\begin{proof}[Proof of the Rationalization Lemma \ref{l:impro}]
By assumption there is a rational function $T\in\Q_\varepsilon(\vec u)$ such that $T(\vec x)\in F(r)-F$.
By the Irreducibility Lemma \ref{l:irre}.a $F(r)=F[r]$.
Hence
$$T(\vec x)=P(r)=p_0+p_1 r+\ldots+p_{k-1} r^{k-1}$$
for some polynomial $P\in F[z]$ of degree less than $k$.
Since $P(r)\not\in F$, there is $l$ such that $0<l<k$ and the coefficient $p_l\in F$ of $z^l$ in $P$ is non-zero.
We have
$$\rho:=p_l r^l=\frac{P(r)+\varepsilon_k^{-l}P(r\varepsilon_k)+\varepsilon_k^{-2l}P(r\varepsilon_k^2)\ldots+
\varepsilon_k^{(1-k)l}P(r\varepsilon_k^{k-1})}k.$$
Define the {\it resolvent polynomial}
$Q(t):=\prod\limits_{\alpha\in S_n}(t-T(\vec x_\alpha))$,
where $S_n$ is the set of all permutations of $\{1,\ldots,n\}$.
Since $T(\vec x)=P(r)$, we have $Q(P(r))=0$.
The coefficients of $Q$ as a polynomial of $t$ are symmetric in $x_1,\ldots,x_n$.
Since $F$ contains elementary symmetric polynomials of $x_1,\ldots,x_n$, it follows that $Q(t)\in F[t]$.
Thus $Q(P(z))\in F[z]$.
Take any $j=1,\ldots,k-1$.
Then by the Conjugation Lemma \ref{l:irre}.b $Q(P(r\varepsilon_k^j))=0$.
Thus $P(r\varepsilon_k^j)=T(\vec x_\alpha)$ for some permutation $\alpha=\alpha_j$.
Hence the above formula for $\rho$ shows that $\rho\in\Q_\varepsilon(\vec x)$.

We have  $\rho^k=p_l^k(r^k)^l\in F$ and $\rho=p_lr^l\in F(r)$.
Since $k$ is a prime and $l$ is not divisible by $k$, there are integers $a$ and $b$ such that $ak+bl=1$.
Since $F(r)$ is a field, we have $r=(r^k)^a(r^l)^b=(r^k)^a\rho^bp_l^{-b}\in F(\rho)$.
Hence $F(r)=F(\rho)$.
\end{proof}

\begin{proof}[Proof of Abel/Ruffini Theorem \ref{t:aberuf}]
We may assume that such that all the degrees of the roots extracted have prime power, i.e.
$F_j=F_{j-1}(r_j)$ for some $r_j\in\C$ such that $r_j^{k_j}\in F_{j-1}$ for some prime $k_j$.
Let us prove the modified statement by  induction on $s$.
Base $s=0$ is obvious.
Let us prove the inductive step.\footnote{This proof is invented by I. Gaiday-Turlov and A. Lvov.
Like everything in this paper, it could have been known earlier.
This proof covers a minor gap in the previous version.}

Take minimal $s$ for which there is a radical extension $\Q_\varepsilon(\vec a)=F_0\subset\ldots\subset F_s\ni x_1$
such that all $k_1,\ldots,k_s$ are prime powers.
Call an integer $j\in\{1,\ldots,s\}$ {\it interesting} if $F_j\cap \Q_\varepsilon(\vec x)\not\subset F_{j-1}$.
By the minimality of $s$, the number $s$ is interesting.
Take the maximal $m\le s$ for which there is a radical extension as above such that $m-1$ is not interesting but $m,m+1,\ldots,s$ are interesting.

If $m=1$, then by the Rationalization Lemma \ref{l:impro} we can consecutively for $j=1,2,\ldots,s$ replace $r_j$ by $\rho_j$  so that
$\Q_\varepsilon(\vec a,\rho_1,\ldots,\rho_j)=F_j$.
For $j=s$ we obtain the required radical extension.

Now assume that $m>1$.
Since $m$ is interesting, by the Rationalization Lemma \ref{l:impro} there is $\rho\in\Q_\varepsilon(\vec x)$ such that $F_m=F_{m-1}(\rho)$.
For the prime $k_m$ we have $\rho^{k_m}\in F_{m-1}$.
Since $m-1$ is not interesting, we have $F_{m-1}\cap \Q_\varepsilon(\vec x)\subset F_{m-2}$.
Hence $\rho^{k_m}\in F_{m-2}$.
Thus $F_{m-2}(\rho)$ is a radical extension of $F_{m-2}$ and $F_{m-2}(\rho,r_{m-1})=F_{m-2}(r_{m-1},\rho)=F_{m-1}(\rho)=F_m$.
So replacement of $r_{m-1},r_m$ by $\rho,r_{m-1}$ gives a radical extension of $F_0$ having the same $s$ and whose each field except $F_{m-2}(\rho)$ coincides with the corresponding field among $F_0,\ldots,F_s$.
For the new extension $m$ is not interesting but $m+1,\ldots,s$ are interesting.
Since $m<s$, this contradicts to the maximality of $m$.
\end{proof}

\comment

\begin{proof}[Attempt of another proof of the analogous result]
We use induction on $s$.
Base $s=0$ is obvious.
Let us prove the inductive step.
Let $\Q_\varepsilon(\vec a)=F_0\subset\ldots\subset F_s$ be a radical extension such that all the roots extracted have prime power.
Apply the inductive step to the radical extension $F_0\subset\ldots\subset F_{s-1}$.
Denote by $Q_0\subset\ldots\subset Q_t$ the obtained radical extension.
If $F_s\cap \Q_\varepsilon(\vec x)= F_{s-1}\cap \Q_\varepsilon(\vec x)$, then this radical extension is as required.
If $F_s\cap \Q_\varepsilon(\vec x)\ne F_{s-1}\cap \Q_\varepsilon(\vec x)$, then take $r$ and prime $k$ such that $F_s=F_{s-1}(r)$.
Take $\rho$ given by the Rationalization Lemma \ref{l:impro}.
Then the radical extension $Q_0\subset\ldots\subset Q_t\subset Q_t(\rho)$ is as reqired???.
Why $Q_t(\rho)=F_s\cap \Q_\varepsilon(\vec x)$?

If $F,A\subset\C$ are fields, $k$ a prime, $\rho\in A$ and $\rho^k\in F$, is $F(\rho)\cap A=(F\cap A)(\rho)$?
\end{proof}

{\it The part of the argument in version 2.}
Assume to the contrary that some root $x_1$ of the equation is contained in some radical extension of
$\Q(\vec a):=\Q(a_{n-1},\ldots,a_0)$.
Then $x_1$ is contained in some radical extension of $\Q_\varepsilon(\vec a)$.
Abbreviate $Q_j:=\Q_\varepsilon(\vec a,r_1,\ldots,r_j)$

Since $s$ is minimal, $x_1\not\in F_{s-1}:=\Q_\varepsilon(\vec a,r_1,\ldots,r_{s-1})$.
Then by the Rationalization Lemma \ref{l:impro} there is $\rho\in\Q_\varepsilon(\vec x)$ such that $\Q_\varepsilon(\vec a,r_1,\ldots,r_s)=F_{s-1}(\rho)$ and $\rho^{k_s}\in F_{s-1}$.

Take a radical extension $\Q_\varepsilon(\vec a,r_1',\ldots,r_{s'}')$ of $\Q_\varepsilon(\vec a)$ containing $\rho$, with minimal $s'$.
Then $s'<s$.
Repeat the previous argument (or proceed by induction on $s$) to obtain a radical extension
$\Q_\varepsilon(\vec a,\rho_1,\ldots,\rho_t)$ of $\Q_\varepsilon(\vec a)$ containing $x_1$ and such that $\rho_j\in \Q_\varepsilon(\vec x)$ for every $j=1,\ldots,t$.
Then $\Q_\varepsilon(\vec a,\rho_1,\ldots,\rho_t)\subset \Q_\varepsilon(\vec x)$ is a radical extension of  $\Q_\varepsilon(\vec a)$ containing $x_1$.
This contradicts to the Ruffini Theorem \ref{t:rufnum}.

\endcomment

\section*{Discussion}

\begin{Remark}[Comparison with other expositions]
The direct expositions of \cite{FT, L, PS, Ko, B, Ka, R, P, St94} were more useful to me (in spite of some drawbacks mentioned below) than `theoretical' expositions in standard textbooks.
The latter start with several hundreds pages of definitions and results whose role in the proof of the insolvability theorem is not clear at the moment of their formulations.
It would not have been possible to write the present note if the above-mentioned direct expositions did not exist.


The above-given proofs of the Ruffini Theorem \ref{t:rufnum} and Theorem \ref{t:galois} are based on \cite[\S9.3]{Ko}, \cite{L, PS} and are similar to \cite{Ay}, \cite[Theorem 5.3]{Bew}\footnote{Cf. Remark \ref{r:ayoub} and footnote \ref{f:form}.
Besides, in \cite{PS} the solvability in radicals of the polynomial $G$ used in the proof of Theorem 4 (p. 219) is not defined.
So instead of Theorem 3 and the first part of the proof of Theorem 4 one needs to use more general results.
These results are not stated.
The exact meaning of `we are dealing with a general polynomial of degree $n$' is not clear (p. 220, after the formula for $\rho_1$).}.
The above-given proofs are different from the cited proofs.


The above-given proof of the Ruffini Theorem \ref{t:rufnum} has the same idea as the proof of \cite{St94}.
This idea is presented above in a more elementary way by looking at the symmetry of a polynomial not at an automorphism of a field.
Still, I like the proof of \cite{St94} and present it (corrected and shortened) in Remark \ref{r:fieaut}.

The exposition here is different from \cite{A, FT, Sk11} (even at the level of formulations, see Remark \ref{r:stat}).
The two approaches are `dual': here the symmetry group of a polynomial is {\it decreased} after the extraction of a radical, while in \cite{A, FT, Sk11}
the group of permutations of roots
is {\it increased} after the extraction of a radical (cf. the proof of Remark \ref{r:fieaut}).


In \cite{Sk17, T} a stronger Kronecker's Theorem is proved; the proof is more complicated.
See also \cite{PC} (so far I did not have time to study that paper).
The exposition of \cite{L} is highly illuminating; it does not claim to be rigorous (and is not).

Before studying the Ruffini Theorem \ref{t:rufnum} a reader might want to learn its simpler {\it real} analogue \cite[\S2.5]{ECG}, \cite[\S8.4.B]{Sk20}.
A real analogue of
Theorem \ref{t:galois} is also simpler \cite[\S6.3]{Sk17}, \cite[\S5.5.3]{ZSS} \cite[\S8.4.E]{Sk20} but uses the ideas of the Irreducibility and Conjugation Lemma \ref{l:irre} in a different way.

See more general remarks and references in \cite[\S2.1]{Sk17}, \cite[\S5.2]{ZSS},
\cite[\S8.1.E and \S8.1.F]{Sk20}.
\end{Remark}


\begin{Remark}[On statements of the insolvability theorems]\label{r:stat}
I was surprised not to find a rigorous statement of the Abel-Ruffini Theorem in Wikipedia
(English, French, German, Italian, Russian, in 2015-2020).
Abel's own paper (see a translation in \cite{P}) does not contain a rigorous formulation in the sense of modern  mathematics (same holds for \cite{P}).
A possible reason is that the Abel-Ruffini Theorem is not so easy to state, see below.
There are different formalizations of the intuitive notion of insolvability in radicals:

$\bullet$ the easiest-to-state insolvability of a specific equation: no root of the equation $x^5-4x+2=0$ is contained in any radical extension of $\Q$
(Galois Theorem),

$\bullet$ an easy-to-state existence of a specific insolvable equation
(Theorem \ref{t:galois}),

$\bullet$ a harder-to-state non-existence of a formula for solution of any equation
(the Abel-Ruffini Theorem, see the Functional and the Formal Abel-Ruffini Theorems of Remarks \ref{r:rufab} and \ref{r:forufab}).


Galois Theorem is much harder than
Theorem \ref{t:galois} and the Abel-Ruffini Theorem.\footnote{`For all that Abel's methods could prove, every particular quintic
equation might be soluble, with a special formula for each equation' \cite{St15}.
This says that even
Theorem \ref{t:galois} is much harder than the Abel-Ruffini Theorem (hence the name of Theorem \ref{t:galois}).
This is unfair from the point of view of modern mathematics because all the ideas of the proof presented here except the existence of algebraically independent numbers were known to Abel or are due to him.
Although proof of the existence uses ideas unfamiliar to Abel, this proof is simple FMPV (=from modern point of view).
So if we name Abel-Ruffini Theorem after Abel and Ruffini in spite of neither having even a precise statement FMPV, it might be fair to name
Theorem \ref{t:galois} `the strong, or the numeric, Abel-Ruffini Theorem' (such a name does not mean that this result was proved by Abel-Ruffini).}
Theorem \ref{t:galois} is easier both to state and to prove than the Functional and the Formal Abel-Ruffini Theorems.
So some readers might want to ignore the statement of the latter.
\end{Remark}





\begin{Remark}[the Functional Abel-Ruffini Theorem]\label{r:rufab}
The exposition of \cite[5.2, 5.3 and Theorem 5.1]{FT}, \cite{A, Sk11} uses the following statement of the Ruffini-Abel Theorem, cf. \cite{Es}.

We may represent solution of a quadratic equation $x^2+px+q=0$ by a sequence of formulas
$$f_1^2=p^2-4q,\quad x=(f_1-p)/2.$$
We may represent solution of a cubic equation $x^3+px+q=0$ by a sequence of formulas
$$f_1^2=\left(\dfrac p3\right)^3+\left(\dfrac q2\right)^2,\quad f_2^3=-\frac q2+f_1,\quad f_3^3=-\frac q2-f_1,\quad x=f_2+f_3.$$
(See comments on the extraction of complex roots in \cite[\S5.2]{FT}.)

These examples motivate the following definition.

Denote the elementary symmetric polynomials in variables $u_1,\ldots,u_n$ by
$$\sigma_1(u_1,\ldots,u_n):=u_1+\ldots+u_n,\quad \ldots,\quad \sigma_n(u_1,\ldots,u_n)=u_1\cdot\ldots\cdot u_n.$$
A {\it radical $n$-formula} is a collection of

$\bullet$ primes $k_1,\ldots,k_s$,


$\bullet$ rational functions $P_0,P_1,\ldots,P_s$ with complex coefficients and in $n,n+1,\ldots,n+s$ variables, respectively, and

$\bullet$ functions $f_1,\ldots,f_s:\C^n\to\C$ (which are not assumed to be continuous)

such that for every numbers $x_1,\ldots,x_n\in\C$
$$\begin{cases} f_1^{k_1}= P_0(\sigma_1,\ldots,\sigma_n) \\
f_2^{k_2}= P_1(\sigma_1,\ldots,\sigma_n,f_1) \\
\ldots\\
f_s^{k_s}= P_{s-1}(\sigma_1,\ldots,\sigma_n,f_1,\ldots,f_{s-1})\\
x_1=P_s(\sigma_1,\ldots,\sigma_n,f_1,\ldots,f_s)
\end{cases}.$$
In these formulas the argument $(x_1,\ldots,x_n)$ of polynomials $\sigma_1,\ldots,\sigma_n$ and functions $f_1,\ldots,f_s$ is omitted;
the equalities are equalities of functions; we assume that the values  $P_j(\sigma_1,\ldots,\sigma_n,f_1,\ldots,f_j)$ are defined.

E.g. there is a radical
2-formula: take
$$s=1,\quad k_1=2,\quad f_1(x_1,x_2)=x_1-x_2,\quad P_0(y_1,y_2)=y_1^2-4y_2\quad\text{and}\quad P_1(y_1,y_2,z_1)=\frac{y_1+z_1}2;$$
$$\text{check that}\quad f_1^2(x_1,x_2)=P_0(x_1+x_2,x_1x_2)\quad\text{and}\quad
x_1=P_1(x_1+x_2, x_1x_2,f_1(x_1,x_2)).$$


{\bf Functional Abel-Ruffini Theorem.} {\it For every $n\ge5$ there is no radical $n$-formula.}


This result follows from
Theorem \ref{t:galois}.
A different proof is presented \cite{A, FT, Sk11}.
That proof uses topological ideas, which could be an advantage for a professional but makes the proof less accessible for a beginner.
(E.g because one has to prove that certain multi-valued functions have the monodromy property \cite[2.10, 2.11]{A}).\footnote{Proof of the Functional Abel Theorem in \cite{Sk15} is incomplete because the definition of an $f$-formula in \cite{Sk15} before Lemma 7 is meaningless for a non-symmetric polynomial $f\in \C[u_1,\ldots,u_n]$.
Merging the fields into $\C(u_1,\ldots,u_n)$ can make the polynomial $y_j^{k_j}-P_{j-1}$ reducible, even if it was irreducible over $\C(\sigma_1,\ldots,\sigma_n)$ for the elementary symmetric polynomials $\sigma_1,\ldots,\sigma_n$ of $u_1,\ldots,u_n$.
E.g. the polynomial $y^2-(u_1+u_2)^2+4u_1u_2$ is reducible over $\C(u_1,u_2)$, although
$y^2-v_1^2+4v_2$ is irreducible over $\C(v_1,v_2)$.
\newline
This gap, as well as the gap mentioned in footnote \ref{f:stil}, are yet another examples that mistakes usually come not from explicitly wrong statements, but from lack of accurate definitions or from use of ill-defined objects.
Cf. \cite{MM}.
For this reason, nowadays accurate definitions are required to recognize a proof as complete.}
\end{Remark}


\begin{Remark}[the Formal Ruffini Theorem]\label{r:ruf}
A {\it rational (or Ruffini) radical $n$-formula} is defined in the same way as a radical $n$-formula of Remark \ref{r:rufab}, except that `functions $f_1,\ldots,f_s:\C^n\to\C$ (which are not assumed to be continuous)' is replaced by `rational functions $f_1,\ldots,f_s\in\C(x_1,\ldots,x_n)$', and `equalities of functions' is replaced by `equalities of rational functions'.

E.g. the previous example of a radical 2-formula is a rational radical 2-formula.


{\bf Formal Ruffini Theorem.} {\it For every $n\ge5$ there is no rational radical $n$-formula.}

This holds by the Ruffini Theorem \ref{t:rufnum}, or can analogously be derived from Lemma \ref{l:radq}.
\end{Remark}

\begin{Remark}[the Formal Abel-Ruffini Theorem]\label{r:forufab}
Let us present a formalization of the statement of the Abel-Ruffini Theorem  in \cite[\S2,\S4]{R}, \cite[Theorem 3]{St94},
\cite[\S1 and Theorem 6.3]{B}.\footnote{\label{f:form} This formalization was not given in \cite{R, St94, B} (so, formally, \cite{R, St94, B} do not contain a rigorous formulation of the Abel-Ruffini Theorem).
Indeed, the objects $\sqrt[p_i]{\alpha_i}$ in \cite[\S2]{R} and $f_{k-1}^{1/m_k}$ in \cite[\S1]{B}, as well as `equals' in the phrase `The adjunction is called {\it radical} if some positive integer power $\alpha^m$ of $\alpha$ equals to an element $f\in F$' \cite[p. 23]{St94} is not defined.
(The object $\alpha_i$ is not a complex number because $E_0$ is $\C(s_1,\ldots,s_n)$ in the simplest formulation of the Abel-Ruffini Theorem there; the object $f_{k-1}$ is an algebraic function not a complex number;
$\alpha^m$ is not equal to $f$ in the only previously defined field $F(\alpha)$ that contains both elements).
The object $E_i^*$ in \cite[\S2]{R} is also not defined, but it might be a typo and $E_i$ is meant.
\newline
In  \cite{St94} it is not specified whether $x_1,\ldots,x_n$ are numbers or variables (since they appear in statements without quantifiers, and since algebraic independent numbers are not mentioned, they have to be variables).
Analogous problem appears in \cite{PS}.
There is no `equation (5.1)'; presumably equation in p. 215 is meant.
It is not defined what is meant by `the coefficients are considered as independent variables over $\Q$'.
It is not defined what is meant by `the roots of a polynomial (5.1) with variable coefficients',
so definition of $\Delta(F)$ is not clear.
Presumably $\Delta=\Q(\sigma_1,\ldots,\sigma_n)$ and $\Delta(F)=\Q(\alpha_1,\ldots,\alpha_n)$ or, in the notation of this note, $\Delta(F)=\Q(u_1,\ldots,u_n)$ for variables $u_1,\ldots,u_n$.
But then it is not defined what is meant by $\sqrt[s_1]{a_1}$, where $a_1\in\Q(u_1,\ldots,u_n)$.
So the main notion of solvability in radicals is not defined (p. 215).
The exact meaning of `we may assume that the roots are independent variables' is not given (p. 220).
\newline
The formalization of $\sqrt[p_i]{\alpha_i}$ suggested in
\cite[definition of a formal radical formula before Lemma 7]{Sk15} in incorrect because
it was not required that the polynomial $y_j^{k_j}-p_{j-1}$ is irreducible over $F_{j-1}(y_j)$.
(That formalization appeared not in the statement but in the proof of the Abel-Ruffini Theorem.)
Such a polynomial is reducible e.g. for $n=s=2$, $p_0(y)=p_1(y)=y^2$.
\newline
For a modern mathematician it is easy, although it does require some accuracy, to formally define $\sqrt[p_i]{\alpha_i}$ (or $f_{k-1}^{1/m_k}$ or `equal' or $\sqrt[s_1]{a_1}$), for some values of $p_i$, $\alpha_i$ not for all values as it seems to be assumed in \cite{R, St94, B, PS}.
(E.g. for $F=\Q(\sqrt[3]2)$ the ring $F[y]/(y^3-2)$ is not a field; in other words, extension of $F$ by $\varepsilon_3\sqrt[3]2$ is not isomorphic to $F[y]/(y^3-2)$.)
See the definition of a formal radical $n$-formula in Remark \ref{r:forufab}.}


Let $n$ be a positive integer.
Denote $F_0:=\C(v_1,\ldots,v_n)$.
{\it A formal radical $n$-formula} is a collection of primes $k_1,\ldots,k_s$ and polynomials
$P_j\in F_0[y_1,\ldots,y_j]$, $j=0,\ldots,s-1$, such that the polynomial $G(y_j):=y_j^{k_j}-P_{j-1}(v_1,\ldots,v_n,y_1,\ldots,y_{j-1})$ is irreducible over $F_{j-1}$, where
$F_1,\ldots,F_s$ are defined inductively by $F_j:=F_{j-1}[y_j]/G(y_j)$.

(Here are details for the definition of $F_j$.
Two polynomials $A,B\in F_{j-1}[y_j]$ are called {\it congruent modulo $G$} if $a-b$ is divisible by $G$.
Let $F_j$ be the set of congruence classes.
Since the polynomial $G$ is irreducible over $F_{j-1}$, the addition and the multiplication on $F_{j-1}$ give an addition and a multiplication on $F_j$ in an obvious way.)

{\it Examples.} (a) For $n=2$ take the formal radical 2-formula
$y_1^2=v_1^2-v_2$ (or, formally, $s=1$, $k_1=2$, $P_0(v_1,v_2)=v_1^2-v_2$).
Denote by $[z]$ the congruence class of $z$.
Then $[v_1+y_1],[v_1-y_1]\in F_1$ are the roots of the equation $x^2-2v_1x+v_2=0$.

(b) For $F_0:=\C(p,q)$  take the formal radical 2-formula $y_1^2=p^3+q^2$, $y_2^3=y_1-q$.
Then
$$\left[y_2-y_2^2(y_1+q)p^{-2}\right],\quad \left[\varepsilon_3y_2-\varepsilon_3^2y_2^2(y_1+q)p^{-2}\right],\quad \left[\varepsilon_3^2y_2-\varepsilon_3y_2^2(y_1+q)p^{-2}\right]\quad \in \quad F_2.$$
are the roots of the equation $x^3+3px+2q=0$
(because $[y_2y_2^2(y_1+q)p^{-2}]=[(y_1-q)(y_1+q)p^{-2}]=[p]\in F_2$).

{\bf Formal Abel-Ruffini Theorem.}
{\it For every $n\ge5$ there is no formal radical $n$-formula for which the equation $x^n-v_1x^{n-1}+\ldots+(-1)^{n-1}v_{n-1}x+(-1)^n v_n=0$ has $n$ distinct roots in $F_s$.}



This follows from
Theorem \ref{t:galois}, or is proved analogously to that result (having a formal radical $n$-formula allows not to consider algebraically independent numbers but to work with variables not with numbers from the start).
This also follows from the Functional Abel-Ruffini Theorem of Remark \ref{r:rufab}.
It would be interesting to know if there is a short converse deduction (i.e. `Formal implies Functional'), i.e. a short proof that
{\it a radical $n$-formula gives a formal radical $n$-formula}.\footnote{This was stated without proof in \cite{Sk15}, which is a gap.
The naive way of constructing a formal radical formula does not work.
Indeed, let $f_1(x_1,x_2)=|x_1-x_2|$ and $f_2(x_1,x_2)=x_1-x_2$.
Then $f_1^2=\sigma_1^2-4\sigma_2$, $f_2^2=f_1^2$ is a radical 2-formula, but a collection $k_1=k_2=2$ and $P_0=v_1^2-4v_2$, $P_1=y_1^2$ is not a formal radical 2-formula. }
\end{Remark}


\begin{Remark}[Alternative proof of a weaker version of
Theorem \ref{t:galois}]\label{r:fieaut}
Here we present a proof from \cite{St94}, filling a gap there and shortening the proof by further stripping the standard approach away.

{\bf Theorem.} {\it For every $n\ge5$ there are $a_0,\ldots,a_{n-1}\in \C$ such that some root of the equation
$x^n+a_{n-1}x^{n-1}+\ldots+a_1x+a_0=0$ is not expressible by radicals from $\{a_0,\ldots,a_{n-1}\}$.}

This is weaker than
Theorem \ref{t:galois}.
It would be interesting to know if this result easily implies Theorem \ref{t:galois}, cf. \cite[Main Lemma (c) (Decomposition) in \S6.4]{Sk17}, \cite[Main Lemma (c) (Decomposition) in \S5.5.4]{ZSS},
\cite[Proof of Lemma 8.4.15.a.]{Sk20}.

For $x_1,\ldots,x_n\in\C$ a field $F\supset\Q(\vec x)$ is called {\it even-symmetric} over a subfield $B\subset F$ if for every even permutation $\alpha$ of $\{1,\ldots,n\}$ there is an automorphism of $F$ which maps every element of $B$ to itself and maps $x_j$ to $x_{\alpha(j)}$ for every $j$.

{\bf Symmetrization Lemma.}
{\it There are $x_1,\ldots,x_n\in\C$ such that for every radical extension $F$ of $\Q_\varepsilon(\vec x)$ there is a radical  extension $\overline F$ of $F$ which is even-symmetric over $\Q_\varepsilon$.}%
\footnote{\label{f:stil} The proof of this lemma given in \cite[Proof of Theorem 1]{St94} is incomplete.
First, the notion of radical expression used in the proof is not defined in \cite{St94}, see footnote \ref{f:form}.
Second, the bijection $\sigma$ in that proof is not defined on the adjoined radical.
If $F(r)\cap\Q_\varepsilon(\vec x)\not\subset F$, we cannot set $\psi_\alpha(r)$ to be any $k$-th power root $r_\alpha$ of $\psi_\alpha(r^k)$ because then $\psi_\alpha$ need not be well-defined, for there could be rational functions $P\in F(z)$ and $Q\in\Q_\varepsilon(\vec u)$ such that $P(r)=Q(\vec x)$ but
$P(r_\alpha)\ne Q(\vec x_\alpha)$.
This problem is resolved above by application of the Rationalization Lemma \ref{l:impro}, and thus by working over $\Q_\varepsilon$ not over $\Q$.
(I am grateful to J. Stillwell for his confirmation that there is an error in \cite{St94}.)
\newline
The Rationalization Lemma \ref{l:impro} represents a significant new idea required to deduce
Theorem \ref{t:galois} from the Ruffini Theorem \ref{t:rufnum}.
However the argument of \cite{St94} can be modified (without introduction of new ideas) to provide a complete proof of the Ruffini Theorem \ref{t:rufnum}.
\newline
The `roots of unity' assumption of \cite[Theorem 2]{St94} is not checked in  \cite[proof of Theorem 3]{St94}.
So the argument in \cite[p. 25 above]{St94} is in fact part of the proof of \cite[Theorem 3]{St94}, although it
is not included in \cite[proof of Theorem 3]{St94}.}

{\it Proof.} Analogously to the proof of the Ruffini Theorem \ref{t:rufnum} take $x_1,\ldots,x_n\in\C$  algebraically independent over $\Q_\varepsilon$.

Then the proof is by the induction on the number $s$ from the definition of a radical  extension.

The base $s=0$ is trivial: for $F=\Q_\varepsilon(\vec x)$ take $\overline F=F$.

In order to prove the inductive step assume that a radical extension $F$ of $\Q_\varepsilon(\vec x)$ is even-symmetric over $\Q_\varepsilon$, that $r\in\C-F$ and $r^k\in F$ for an integer $k$.
For every even permutation $\alpha$ take an automorphism $\psi_\alpha:F\to F$ which maps every element of $\Q_\varepsilon$ to itself and $x_j$ to $x_{\alpha(j)}$ for every $j$.

If $F(r)\cap\Q_\varepsilon(\vec x)\not\subset F$, then by the Rationalization Lemma \ref{l:impro} there is $\rho\in\Q_\varepsilon(\vec x)$ such that $\rho^k\in F$ and $F(\rho)=F(r)$.
Take a rational function $P\in \Q_\varepsilon(\vec u)$ such that $\rho=P(\vec x)$.
Extend $\psi_\alpha$ to $\overline F:=F(\rho)=F(r)$ by setting $\psi_\alpha(\rho):=P(\vec x_\alpha)$.
This is well-defined by the Irreducibility Lemma \ref{l:irre}.a because
$\psi_\alpha(P(\vec x)^k)=P(\vec x_\alpha)^k$.

If $F(r)\cap\Q_\varepsilon(\vec x)\subset F$, then for every even permutation $\alpha$ take any $r_\alpha\in\C$ such that $r_\alpha^k=\psi_\alpha(r^k)$.
Take a minimal set $\beta_1,\ldots,\beta_t$ of even permutations such that
$\overline F:=F(r_{\beta_1},\ldots,r_{\beta_t})$ contains $r_\alpha$ for every even permutation $\alpha$.
Extend $\psi_\alpha$ to $\overline F$ by setting $\psi_\alpha(r_{\beta_j}):=r_{\alpha\beta_j}$.
The extension is well-defined by the Irreducibility Lemma \ref{l:irre}.a.

The extended $\psi_\alpha$ maps every element of $\Q_\varepsilon$ to itself and maps $x_j$ to $x_{\alpha(j)}$ for every $j$.
The inductive step is proved.
\qed


{\bf Extraction of Radical Lemma} (analogue of Lemma \ref{l:radq}).
{\it Assume that $n\ge5$ and $k$ are integers, $x_1,\ldots,x_n\in\C$, a field $F\supset\Q(\vec x)$ is even-symmetric over a subfield $B\subset F$, \ $r\in F-B$ and $r^k,\varepsilon_k\in B$.%
\footnote{Here we do not assume that are $x_1,\ldots,x_n$ algebraically independent over $\Q_\varepsilon$, although we apply the Lemma in that situation.}
Then $F$ is even-symmetric over $B(r)$.}

{\it Proof.}
Let $a,b,c,d,e$ be arbitrary different elements of $\{1,\ldots,n\}$.
It suffices to prove that $F$ is $(abc)$-symmetric over $B(r)$, i.e. that there is an automorphism of $F$ which maps every element of $B(r)$ to itself and $x_j$ to $x_{(abc)(j)}$ for every $j$.
Recall that
$$(abc)=(dac)^{-1}(ceb)^{-1}(dac)(ceb).$$
Since $F$ is even-symmetric over $B$, there are  automorphisms $\alpha,\beta$ of $F$ such that

$\bullet$ $\alpha(x_d,x_a,x_c)=(x_a,x_c,x_d)$,

$\bullet$ $\beta(x_c,x_e,x_b)=(x_e,x_b,x_c)$,

$\bullet$ $\alpha(y)=y$ for every $y\in B$ or $y=x_j$, $j\not\in\{d,a,c\}$, and

$\bullet$ $\beta(y)=y$ for every $y\in B$ or $y=x_j$, $j\not\in\{c,e,b\}$.

Since $r^k\in B$, we have $\alpha(r^k)=\beta(r^k)=r^k$.
Hence there are integers $p,q$ such that $\alpha(r)=\varepsilon_k^pr$ and $\beta(r)=\varepsilon_k^qr$.
Since $\varepsilon_k\in B$, we have $\alpha^{-1}(r)=\varepsilon_k^{-p}r$ and $\beta^{-1}(r)=\varepsilon_k^{-q}r$.
Therefore $\alpha^{-1}\beta^{-1}\alpha\beta(r)=r$.
Hence the automorphism $\alpha^{-1}\beta^{-1}\alpha\beta$ of $F$ is as required.
\qed

{\it Alternative proof of the above Theorem} (cf. proof of the Ruffini Theorem \ref{t:rufnum}).
Take $x_1,\ldots,x_n\in\C$ given by the Symmetrization Lemma.
Denote the coefficients of the unitary polynomial with roots $x_1,\ldots,x_n$ by
$a_{n-1},\ldots,a_0$.
Assume to the contrary that $x_1,\ldots,x_n$ are contained in some radical extension of
$\Q(\vec a):=\Q(a_{n-1},\ldots,a_0)$.
Then $x_1,\ldots,x_n$ are contained in some radical extension $F$ of $\Q_\varepsilon(\vec a)$.
Hence $F$ is also a radical extension of $\Q_\varepsilon(\vec x)$.
Therefore there is a radical extension $\overline  F$ of $F$ which is even-symmetric over $\Q_\varepsilon$.
Since $a_{n-1},\ldots,a_0$ are symmetric polynomials of $x_1,\ldots,x_n$, the field $\overline  F$ is  even-symmetric over $\Q_\varepsilon(\vec a)$.
Since $F$ is a radical extension of $\Q_\varepsilon(\vec a)$, the field $\overline  F$ is also a radical extension of $\Q_\varepsilon(\vec a)$, i.e.  $\overline  F=\Q_\varepsilon(\vec a,r_1,\ldots,r_s)$ for some $r_1,\ldots,r_s$.
By induction on $j$ using Extraction of Radical Lemma
one proves that $\overline F$ is even-symmetric over $\Q_\varepsilon(\vec a,r_1,\ldots,r_j)$ for every $j=0,1,\ldots,s$.
For $j=s$ this is a contradiction because even permutation $(13)(12)=(123)$ carries $x_1$ to $x_2\ne x_1$.
\end{Remark}

\begin{Remark}[A relation of the proofs in this paper and in \cite{Ay}]\label{r:ayoub}
The proofs are very close.
The exposition of \cite{Ay} is more general.
More precisely, the Rationalization Lemma \ref{l:impro} corresponds to \cite[Lemma in p. 399]{Ay};
the objects $\Q_\varepsilon(\sigma_1,\ldots,\sigma_n),\Q_\varepsilon(\vec x),F,F(r),r,\rho$ correspond to $H,F,M,L,\lambda,\mu$.

The following non-trivial points (a,b) of the proof are missing in \cite{Ay}.
(Here and below `non-trivial' means `non-trivial for a non-specialist'; these points might be clear for a specialist.)

(a) The existence of automorphisms $\sigma_k$ in \cite[p. 399, second paragraph of the proof of Lemma]{Ay}
is non-trivial (as opposed to the uniqueness of $\sigma_k$; cf. footnote \ref{f:stil}).
This existence follows by the Conjugation Lemma \ref{l:irre}.b not even stated in \cite{Ay}.

(b) The normality of $F$ over $H$ is not defined in \cite{Ay} (see \url{https://en.wikipedia.org/wiki/Normal_extension}).
In the application \cite[Lemma in p. 399]{Ay} checking the normality is omitted and is non-trivial.
(The use of the normality of $F$ over $H$ in \cite[proof of Lemma in p. 399]{Ay} corresponds to the use of Galois resolvent $Q(t)$ in the proof of the Rationalization Lemma \ref{l:impro}.)



(c) Certain systems of linear equations are easier to solve explicitly than by reference to their determinants, see the explicit formula for $\rho$ in the proof of the Rationalization Lemma \ref{l:impro}.
\end{Remark}

\begin{Remark}[A solvability criterion and an algorithm]\label{r:alg}
{\it There is an algorithm deciding, for $a_{n-1},\ldots,a_0\in\Q$, whether all the roots of the equation $x^n+a_{n-1}x^{n-1}+\ldots+a_1x+a_0=0$
are expressible by radicals from $\{1\}$.}


This is implied by the following criterion, together with an estimation on the number of operations.
(This estimation can easily be extracted from the proof,
the idea is to observe that
the `symmetry subgroup' of $S_n$
cannot be changed more than $\log_2n!<n\log_2n$ times.)


{\bf Galois Solvability Criterion (conjecture).}
{\it For every $a_{n-1},\ldots,a_0\in\Q$ all the roots of the equation $A(x):=x^n+a_{n-1}x^{n-1}+\ldots+a_1x+a_0=0$ are expressible by radicals from $\{1\}$
if and only if a set of degree 1 polynomials over $\Q$
can be obtained from $\{A\}$ using the following operations:

$\bullet$ (factorization) if one of our polynomials equals to $P_1P_2$ for some non-constant $P_1,P_2\in\Q[x]$, then replace $P_1P_2$ by $P_1$ and $P_2$;

$\bullet$ (extracting a root) if one of our polynomials equals to $P(x^k)$ for some $P\in\Q[x]$, then replace $P(x^k)$ by $P(x)$;

$\bullet$ (taking Galois resolvent) replace one of our polynomials $P$ by the polynomial
$$\prod\limits_{\alpha\in\Sigma_q} (x-\varepsilon_qy_{\alpha(1)}-\varepsilon_q^2y_{\alpha(2)}-\ldots-\varepsilon_q^ky_{\alpha(q)}),$$
where $y_1,\ldots,y_q$ are all the roots of $P$.
(The coefficients of this product are symmetric in $y_1,\ldots,y_k$, so they are rational, i.e. $y_1,\ldots,y_k$ are `not required' to calculate the coefficients.)}

I would be grateful if a specialist in algebra could confirm that this criterion is correct (and is equivalent to the Galois Solvability Criterion in its usual textbook formulation, please give a reference), or describe required changes.
(I asked some specialists since July 2017, but so far obtained no answer.)

I conjecture that for every $a_{n-1},\ldots,a_0\in\C$ analogous results holds for $a_{n-1},\ldots,a_0\in\C$ with $\{1\}$ replaced by   $\{1,a_{n-1},\ldots,a_0\}$.
\end{Remark}

\end{document}